\documentclass[12pt]{amsart}
\usepackage{amsbsy}
\usepackage{graphicx,epsfig,subfigure,psfrag}
\begin{document}

\newtheorem{theorem}{Theorem}
\newtheorem{proposition}{Proposition}
\newtheorem{lemma}{Lemma}
\newtheorem{corollary}{Corollary}
\newtheorem{definition}{Definition}
\newtheorem{remark}{Remark}
\newcommand{\tex}{\textstyle}
\numberwithin{equation}{section} \numberwithin{theorem}{section}
\numberwithin{proposition}{section} \numberwithin{lemma}{section}
\numberwithin{corollary}{section}
\numberwithin{definition}{section} \numberwithin{remark}{section}
\newcommand{\ren}{\mathbb{R}^N}
\newcommand{\re}{\mathbb{R}}
\newcommand{\n}{\nabla}
\newcommand{\iy}{\infty}
\newcommand{\pa}{\partial}
\newcommand{\fp}{\noindent}
\newcommand{\ms}{\medskip\vskip-.1cm}
\newcommand{\mpb}{\medskip}
\newcommand{\AAA}{{\bf A}}
\newcommand{\BB}{{\bf B}}
\newcommand{\CC}{{\bf C}}
\newcommand{\DD}{{\bf D}}
\newcommand{\EE}{{\bf E}}
\newcommand{\FF}{{\bf F}}
\newcommand{\GG}{{\bf G}}
\newcommand{\oo}{{\mathbf \omega}}
\newcommand{\Am}{{\bf A}_{2m}}
\newcommand{\CCC}{{\mathbf  C}}
\newcommand{\II}{{\mathrm{Im}}\,}
\newcommand{\RR}{{\mathrm{Re}}\,}
\newcommand{\eee}{{\mathrm  e}}
\newcommand{\LL}{L^2_\rho(\ren)}
\newcommand{\LLL}{L^2_{\rho^*}(\ren)}
\renewcommand{\a}{\alpha}
\renewcommand{\b}{\beta}
\newcommand{\g}{\gamma}
\newcommand{\G}{\Gamma}
\renewcommand{\d}{\delta}
\newcommand{\D}{\Delta}
\newcommand{\e}{\varepsilon}
\newcommand{\var}{\varphi}
\newcommand{\lll}{\l}
\renewcommand{\l}{\lambda}
\renewcommand{\o}{\omega}
\renewcommand{\O}{\Omega}
\newcommand{\s}{\sigma}
\renewcommand{\t}{\tau}
\renewcommand{\th}{\theta}
\newcommand{\z}{\zeta}
\newcommand{\wx}{\widetilde x}
\newcommand{\wt}{\widetilde t}
\newcommand{\noi}{\noindent}
\newcommand{\uu}{{\bf u}}
\newcommand{\xx}{{\bf x}}
\newcommand{\yy}{{\bf y}}
\newcommand{\zz}{{\bf z}}
\newcommand{\aaa}{{\bf a}}
\newcommand{\cc}{{\bf c}}
\newcommand{\jj}{{\bf j}}
\newcommand{\ggg}{{\bf g}}
\newcommand{\UU}{{\bf U}}
\newcommand{\YY}{{\bf Y}}
\newcommand{\HH}{{\bf H}}
\newcommand{\GGG}{{\bf G}}
\newcommand{\VV}{{\bf V}}
\newcommand{\ww}{{\bf w}}
\newcommand{\vv}{{\bf v}}
\newcommand{\hh}{{\bf h}}
\newcommand{\di}{{\rm div}\,}
\newcommand{\ii}{{\rm i}\,}
\newcommand{\inA}{\quad \mbox{in} \quad \ren \times \re_+}
\newcommand{\inB}{\quad \mbox{in} \quad}
\newcommand{\inC}{\quad \mbox{in} \quad \re \times \re_+}
\newcommand{\inD}{\quad \mbox{in} \quad \re}
\newcommand{\forA}{\quad \mbox{for} \quad}
\newcommand{\whereA}{,\quad \mbox{where} \quad}
\newcommand{\asA}{\quad \mbox{as} \quad}
\newcommand{\andA}{\quad \mbox{and} \quad}
\newcommand{\withA}{,\quad \mbox{with} \quad}
\newcommand{\orA}{,\quad \mbox{or} \quad}
\newcommand{\atA}{\quad \mbox{at} \quad}
\newcommand{\onA}{\quad \mbox{on} \quad}
\newcommand{\ef}{\eqref}
\newcommand{\ssk}{\smallskip}
\newcommand{\LongA}{\quad \Longrightarrow \quad}
\def\com#1{\fbox{\parbox{6in}{\texttt{#1}}}}
\def\N{{\mathbb N}}
\def\A{{\cal A}}
\newcommand{\de}{\,d}
\newcommand{\eps}{\varepsilon}
\newcommand{\be}{\begin{equation}}
\newcommand{\ee}{\end{equation}}
\newcommand{\spt}{{\mbox spt}}
\newcommand{\ind}{{\mbox ind}}
\newcommand{\supp}{{\mbox supp}}
\newcommand{\dip}{\displaystyle}
\newcommand{\prt}{\partial}
\renewcommand{\theequation}{\thesection.\arabic{equation}}
\renewcommand{\baselinestretch}{1.1}
\newcommand{\Dm}{(-\D)^m}

\title
{\bf On  blow-up shock  waves for a nonlinear\\
 PDE
 associated with Euler
equations}

\author{
V.A.~Galaktionov}

\address{Department of Mathematical Sciences, University of Bath,
 Bath BA2 7AY, UK}
\email{vag@maths.bath.ac.uk}



  \keywords{2D Euler equations, hyperbolic system, shock and rarefaction  waves,
similarity solutions with blow-up swirl.
  {\bf Submitted to:} Math. Theory of Fluid Dyn.}
 \date{\today}




 \begin{abstract}

The following nonlinear PDE:
 $$
 \tex{
  \big(\frac{u_t}u\big)_t= \frac{u_x u_t}u \quad \mbox{in\,\, $\re \times
  \re_+$},
  }
  $$
  derived from 2D Euler equations, is shown to
 admit  smooth similarity solutions,
  which as $t \to T^-$ create  shocks of the type
 ${\rm sign}\, x$ and $H(-x)$ and also other blow-up singularities.
  Some of the blow-up solutions are shown to admit unique
  extension beyond, for $t>T$.
 All the similarity reductions lead to singular boundary value problems
 for 2nd-order ODEs. More complicated solutions with blow-up
 angular swirl are discussed.

\end{abstract}

\maketitle


\section{Introduction: Euler equations and a related PDE with shock waves}
 \label{S1}


\subsection{Euler equations and shocks: beginning of derivation}

Consider the 2D {\em Euler equations} (the EEs) of inviscid
incompressible fluids
  \be
  \label{NS1}
 \uu_t +(\uu \cdot \n)\uu=- \n p, \quad \di \uu=0 \inB \re^2 \times
 (0,T) \quad(\uu=(u_1,u_2)^T),
  \ee
 with   bounded initial
divergence-free  $\uu_0$
 of {\em infinite} kinetic $L^2$-energy.
 Using Leray's formulation with the projector ${\mathbb P}=I- \n
 \D^{-1} (\n \cdot)$, \ef{NS1} is written in the form
 \cite[p.~30]{Maj02}
  \be
  \label{EE1}
 \begin{matrix}
  \uu_t + (D \uu)\,\uu
   = \BB (\uu)  \equiv \frac 1{2 \pi} \,
  \int\limits_{\re^2} \frac {\xx-\yy}{|\xx-\yy|^2} \, {\rm tr}\, (\n
  \uu(\yy,t))^2 \, {\mathrm d}\yy \whereA
   \\
   {\rm tr}\, (\n
  \uu)^2  = \sum\limits_{(i,j)} (u_j)_{x_i}(u_i)_{x_j} \andA
  D \uu = \left[
   \begin{matrix}
 (u_1)_{x_1} \,\,\, (u_1)_{x_2} \ssk\\
 (u_2)_{x_1} \,\,\,(u_2)_{x_2}
  \end{matrix}
  \right]
 \end{matrix}
   \ee
   is  the  $2\times 2$ Jacobian matrix of the solenoidal vector field
$\uu$.

 It is well-known that,
 even in dimension 2 (to say nothing about the 3D ones),
 the standard notion of weak solutions is not sufficient for
 establishing the unique solvability of the problem in the case of the unbounded initial
  vorticity, where
 ${\rm curl}\, \uu_0 \not \in L^\iy$. See Majda--Bertozzi \cite[\S~8]{Maj02} and recent
 surveys in Bardos--Titi
 \cite{Bar07}, Constantin  \cite{Con07}, and Ohkitani \cite{Ohk08}
  for further details and key references, as well as
 Pomeau {\em et al} \cite{Pom08} for discussion of somehow
  related open questions on shock waves in hyperbolic systems.
 It seems obvious that
 natural reasons for this difficulty are two fold. Firstly,
 there is no still a proper definition of ``entropy solutions" for
 the EEs admitting sufficiently strong ``shock  waves". In other words,
  general principles of
formation of various singularities of shock wave types, as well as
$L^\iy$-blow-up, are in general unknown. Secondly, on the other
hand, even the  corresponding system without the nonlocal term
 \be
 \label{HS1}
\uu_t + (D \uu) \uu
   = 0 \inB \re^2 \times \re_+
    \ee
 is a rather difficult divergent hyperbolic system $2 \times 2$ with two
independent spatial variables $\xx=(x_1,x_2)$, so it does not
automatically fall into the scope of modern entropy theory
regardless of its definite recent achievements; see  Bressan
\cite{Bres} and Dafermos \cite{Daf}. Of course, some classes of
solutions of systems such as \ef{HS1} can be represented {\em via}
 characteristics (see interesting  physically motivated
formal discussions in \cite{Pom08}), though, in general,  choosing
right weak ``entropy" solutions is not straightforward  and  is an
open problem.

\ssk

 Our goal is to introduce a new PDE model, which may reflect some features of formation
of finite-time singularities in systems such as \ef{NS1} (and
partially \ef{HS1}) and, in addition, having vorticity-swirl
structure that can be somehow adequate to the  EEs. To this end,
to discuss possible ways of shock wave formation for the EEs, let
us write down \ef{NS1} in the polar coordinates $(r, \var)$,
 \be
 \label{1}
  \left\{
   \begin{matrix}
   U_t+UU_r + \frac 1r\, V U_\var- \frac 1r\, V^2=-P_r,
  \quad\ssk\\
   V_t + UV_r + \frac 1r\, VV_\var + \frac 1r\, UV=- \frac 1r \,
   P_\var,
    \ssk\\
   U_r+ \frac 1r\, U + \frac 1r\, V_\var=0.
   \qquad\qquad\qquad\,\,\,
   \end{matrix}
   \right.
   \ee
 Though \ef{1} is a system of three PDEs, which is equivalent to the nonlocal equation \ef{EE1},
  we consider a
 rather hypothetical situation of
 formation of  shocks
 created by local differential operators only\footnote{Obviously, this is
 not the case for 3D unbounded $L^\iy$-singularities, where quadratic
 nonlocal terms must be involved and this settles the core of such
 a remarkable open problem, \cite{Con07}.}, though, as we have
 seen, this scenario is expected to be also relevant to the
 hyperbolic system \ef{HS1}.

\ssk

 Let us then more precisely define a manifold of solutions
 with shocks to be studied:

\ssk

 (i) Formation of {\em bounded} shocks (or bounded {\em rarefaction waves}) is not essentially affected by the pressure terms,
 i.e., by the nonlocal members (however, the div-free
 equation remains important); and

\ssk

 (ii) As in the classic case \ef{3} below,
two  quadratic first-order operators   in the first
 and in the second  equations in \ef{1} are leading together with two differential
 terms in the div-free equation (the latter is reasonable for finite shocks):
 \be
 \label{2}
  \left\{
   \begin{matrix}
   U_t+UU_r +...=0, \,\,\,\ssk\\
   V_t +  \frac 1r\, VV_\var +...=0,  \ssk\\
   U_r+ \frac 1r\, V_\var +...=0.\,\,\,
   \end{matrix}
   \right.
   \ee
Thus,  we omit other operators that are expected to be bounded on
such singularities.

\ssk

 Replacing $V$ by $r V$ for convenience yields
 \be
 \label{21}
 V \mapsto r V \LongA
  \left\{
   \begin{matrix}
   U_t+UU_r +...=0, \,\,\,\ssk\\
   V_t +  VV_\var +...=0, \,\, \ssk\\
   U_r+ V_\var +...=0.\,\,\,\,\,
   \end{matrix}
   \right.
   \ee

The first two equations are indeed the famous 1D {\em Euler
equations} from gas dynamics,
 \be
 \label{3}
 u_t + u u_x=0 \inB \re \times \re_+, \quad u(x,0)=u_0(x) \in L^1
 \cap L^\iy,
  \ee
  whose entropy theory  was created by Oleinik \cite{Ol1, Ol59} and Kruzhkov
\cite{Kru2} (equations in $\ren$) in the 1950--60s; see details on
the history, main results, and modern developments in the
well-known monographs \cite{Bres, Daf, Sm}\footnote{First study of
finite time singularities as shock waves of quasilinear equations
  was performed  by Riemann in 1858 \cite{Ri58} (Riemann's method and invariants are originated therein); see
\cite{Chr07, Pom08} for details. The implicit solution of the
problem \ef{3}, $u(x,t)=u_0(x-u(x,t)t)$ (containing the key
 wave ``overturning" effect) was obtained earlier by Poisson in 1808
\cite{Poi08}; see \cite{Pom08}.}. Typical formation of shocks for
\ef{3} as $t \to 0^-$ (the blow-up time here is $T=0$) is
described by the canonical self-similar blow-up solution
 \be
 \label{4}
  u_{-}(x,t)= \left\{
 \begin{matrix}
  1 \quad \mbox{for} \quad x<t,\\
 \,\,\,\,\frac xt \quad \mbox{for} \quad |\frac xt|<1,\\
 -1 \quad \mbox{for} \quad x > -t.
  \end{matrix}
   \right.
   \ee
Therefore, the following steady shock is formed:
 \be
 \label{5}
  \tex{
   \lim_{t \to 0^-}\, u_{-}(x,t) = S_-(x) \equiv -{\rm sign}\,
   x,
   }
   \ee
 with convergence in $L^1(\re)$ and uniformly on any subset bounded away from the origin.
Moreover, we have that, as $t \to 0^-$,
 \be
 \label{6}
  \tex{
  D_x u_{-}(x,t)= \frac 1t \to - \infty \,\,\, \mbox{for
  any}\,\,\, |x| <|t|
  \andA D_x u_{-}(x,t) \to - 2 \d(x).
  }
  \ee

Indeed, such a standard scenario of shock wave formation is not
applicable for the system in \ef{21}, since the last div-free
equation
 \be
 \label{7}
U_r+ V_\var +...=0
 \ee
prohibits simultaneous blow-up of derivatives  $U_r$ and $V_\var$
to $-\iy$. Therefore, the system \ef{21} needs another more
involved treatment of possible (if any) shock waves. We next
introduce such a nonlinear PDE model.

\subsection{A simplified PDE model and main results}

Thus, we consider  first two equations in \ef{21} on the
``critical manifold" \ef{7}, i.e.,
\be
 \label{8}
  \left\{
   \begin{matrix}
   U_t+UU_r=0, \ssk\\
   V_t +  VV_\var=0\,
   \end{matrix}
 \right.
   \Big|_{U_r+V_\var=0}
    \LongA
    \fbox{$
\left\{
   \begin{matrix}
   U_t=UV_\var, \ssk\\
   V_t =  VU_r,
   \end{matrix}
 \right.
 $}
   \ee
 with given bounded smooth initial data $(U_0,V_0)(r,\var)$.
Writing this in the form
\be
 \label{9}
  \left\{
   \begin{matrix}
   (\ln |U|)_t= V_\var, \ssk\\
   (\ln |V|)_t= U_r,\,
   \end{matrix}
 \right.
 \ee
we find on integration of the second equation that
 \be
 \label{10}
  \tex{
\ln \big|\frac {V(t)}{V_0} \big|= \int_0^t U_r(s)
 \LongA V(t)= V_0 {\mathrm e}^{\int_0^t U_r(s)}
.
 }
  \ee
Substituting this $V(t)$ into the first equation for $U$ yields
the following nonlocal PDE:
 \be
 \label{11}
  \tex{
  U_t= U V_{0\var}{\mathrm e}^{\int_0^t U_r(s)} +
   U V_0 \, \big({\mathrm e}^{\int_0^t U_r(s)}\big)\, {\int_0^t U_{r\var}(s)}
   \equiv J_1+J_2.
    }
    \ee
 According to \ef{11}, there exist at least two different mechanisms of
 shock wave formation corresponding to nonlinear integral
 operators $J_1$ and $J_2$, and possibly, other types via their
 nonlinear interactions. Note that the problem corresponding to
 $J_2$ is assumed to be truly 2D (and hence difficult), while that for $J_1$ can be
 reduced to 1D.

 We begin with the simplest scenario corresponding to $J_1$, for which
 we choose initial data $V_0$ such that, in a small neighbourhood
 of a possible shock at $(0,0)$,
  \be
  \label{12}
   V_{0\var} \approx +1 \andA V_0 \approx 0.
    \ee
    This locally yields the equation
  \be
  \label{13}
  \tex{
 U_t= U {\mathrm e}^{\int_0^t U_r(s)}
 \quad \mbox{or} \quad \ln \big( \frac {U_t}U \big)=
{\int_0^t U_r(s)}, \quad \mbox{where we assume} \,\,\,
\frac{U_t}U>0.
 }
 \ee
 Differentiating this equation in $t$ and, for convenience,  renaming the variables
  \be
  \label{14}
 r \mapsto x \andA U(r,t) \mapsto u(x,t)
  \ee
  yields the necessary equation for such shock wave and blowing up
  formation:
   \be
   \label{15}
    \tex{
    \big( \frac {u_t}u\big)_t= \frac { u_x \,u_t}u
    \quad \mbox{or}
    \quad
    \fbox{$
     u u_{tt}-(u_t)^2=u u_x u_t.
     $}
     }
     \ee
Written in the form with the evolution time-variable $x$ and $t$
being the space variable,
 \be
 \label{151}
  \tex{
  u_x= \frac 1{u_t} \big[ u_{tt}- \frac{(u_t)^2}u\big],
   }
   \ee
 the equation becomes
  \be
  \label{152}
   \mbox{of parabolic type in the domain of time-monotonicity
   $\{u_t>0\}$}
   \ee
   and parabolic backward in time otherwise.
 In what follows, we will need carefully control the parabolicity
 condition \ef{152} in order to guarantee the well-posedness of
 our shock and blow-up phenomena.
    In the original
  independent  variables $\{x,t\}$, which assume two initial
  conditions at $t=0$, \ef{15} exhibits different evolution properties.
Note that, in terms of $t$ as the time variable, this equation
 \be
 \label{151N}
  \tex{
  u_{tt}= \frac {u u_x u_t +(u_t)^2}u
   }
   \ee
 is
in the normal form in the non-singular set $\{u \not = 0\}$,
 so it
obeys there the Cauchy--Kovalevskaya Theorem \cite[p.~387]{Tay}.
Hence, for analytic initial data $u(x,0) \not =0$ and $u_t(x,0)$,
there exists a unique local in time analytic solution $u(x,t)$. On
the other hand, any existing solution $u(x,t)$, which is analytic
in $\{u \not = 0\}$, has a unique local analytic continuation at
any point, where $u \not = 0$.

Let us make an important comment concerning the type of this new
PDE. Written as a first-order system for $U=(u,v)^T$, with
$v=u_t$, \ef{15} takes the form
 \be
 \label{HS10}
  \tex{
  U_t= A(U) U_x +F(U) \whereA
  A(U) =\left[
  \begin{matrix}
  0\,\,\,0\\
  v \,\,\,0
  \end{matrix}
  \right] \andA
  F(U)= \left[
  \begin{matrix}
  v\\
   \frac{v^2}u
  \end{matrix}
  \right].
   }
   \ee
It follows that the characteristic equation for the matrix $A(U)$
has the trivial form
 \be
 \label{ch1}
 {\rm det}\,(A(U)-\l I_2)=0 \LongA \l^2=0.
  \ee
Therefore the eigenvalues always coincide, $\l_\pm=0$,
 so that
 \be
 \label{ch2}
 \mbox{(\ref{HS10}) is not  a strictly hyperbolic system},
  \ee
 and indeed it is {\em uniformly degenerate}.
 In view of
\ef{ch2},
 application of classic theory of
hyperbolic systems (see \cite{Bres, Daf}), which recently got a
fundamental progress, is not possible.
 In any case, it is worth recalling  that
formation of singularities for \ef{HS10} or \ef{15} and possible
extensions of blow-up solutions beyond singularities remain the
key (and absolutely  unavoidable for any difficult nonlinear PDE
model) question to be addressed in what follows.


Regardless such nice local regularity and analyticity properties
of the PDE
 \ef{15},
 our main goal is to show that \ef{15} admits rather
standard ``quasi-stationary" self-similar shock waves (Section
\ref{S2}) and other blow-up singularities (Section \ref{S2N}). In
Section \ref{S31}, following the strategy in \cite{GalndeIII}, we
develop a local theory of formation of shock waves from continuous
data and show that the PDE \ef{15} admits a unique similarity
continuation after blow-up. In other words, this implies that
 uniqueness is not lost locally at the
singularity blow-up points, and this links \ef{15} with the
classic model  \ef{3}. However, a proper definition of entropy
solutions to create a global existence-uniqueness theory for
\ef{15} remains a difficult open problem (it is not clear if such
one can be developed in principle).

\subsection{On blow-up similarity solutions with swirl}

Let us for a while return to the full boxed system in \ef{8},
which is a difficult object, so we truly and  desperately need a
simplified model. For instance, in general, \ef{8} admits
complicated singular solutions with {\em blow-up angular swirl}.
For the NSEs and EEs, such a mechanism was proposed in
\cite{GalJMP}, so we omit details, and prescribe those similarity
solutions as follows: for fixed  $\a \in (0,1)$ and $\s \in \re$,
 \be
 \label{sw1}
 \begin{matrix}
 U(r,\var,t)=(-t)^{-\a}f(y,\mu), \quad
 V(r,\var,t)=(-t)^{-1}g(y,\mu), \ssk\ssk\\
y= \frac r{(-t)^\b}, \quad \b=1-\a, \quad
 \mu=\var + \s \ln(-t), \qquad\qquad\quad
 \end{matrix}
 \ee
where we introduce blow-up rotation in the angular direction being
a {\em logarithmic travelling wave}. Substituting \ef{sw1} into
\ef{8} yields the following first-order system:
 \be
 \label{sw2}
  \left\{
  \begin{matrix}
  \a f+ \b f'_y y - \s f'_\mu= fg'_\mu,\ssk\\
  g+ \b g'_y y - \s g'_\mu= g f'_y.\,\,\,
   \end{matrix}
   \right.
   \ee
This is a difficult system, which we do not intent to study it
here, though indeed some its features will be detected for a
reduced simpler PDE model in \ef{15}.

\section{Shock waves by blow-up similarity solutions}
 \label{S2}

\subsection{Existence of self-similar shock formation}

Similar to \ef{4}, we look for blow-up similarity solutions of
\ef{15} of the form
 \be
 \label{16}
  \tex{
  u_{-}(x,t)= f(y), \quad y= \frac x{(-t)} \whereA
  \big( \frac {f'y}y\big)'+ \frac{f'}f= \frac{(f')^2}f.
   }
   \ee
The ODE in \ef{16} can be easily reduced to the following one:
 \be
 \label{17}
  y f f''=(f')^2(y+f)-2 ff', \quad f(\pm \iy)=\pm 1,
   \ee
   which is degenerate at $y=0$ and at any zero $f=0$ of
   solutions. Note that, unlike \ef{4}, the conditions at infinity
   have been changed, so that \ef{16} is supposed to deliver the
   opposite shock:
 \be
 \label{18}
  \tex{
   \lim_{t \to 0^-}\, u_{-}(x,t) = S_+(x) \equiv {\rm sign}\,
   x.
   }
   \ee
In view of the obvious symmetry of \ef{17}:
 \be
 \label{19}
  \left\{
   \begin{matrix}
    y \mapsto -y,\\
    f \mapsto -f
 \end{matrix}
 \right.
  \ee
  it suffices to solve the following problem:
\be
 \label{20}
  y f f''=(f')^2(y+f)-2 ff' \forA y>0; \quad f(0)=0, \,\,\, f(+\iy)=
  1.
   \ee
Recall that according to \ef{13} we have to check the sign of
 \be
 \label{ss1}
  \tex{
  \frac{u_t}u \equiv (-t)^{-1} \frac{f' y}f>0
  }
  \ee
  which is true if $f(y)>0$ and  $f'(y)>0$ for all $y>0$.
 Here, \ef{20} is a pretty standard second-order ODE, so we
 state the final result:

 \begin{theorem}
 \label{T1}
 The problem $\ef{20}$ admits a unique $C^2([0, \iy))$ strictly monotone
 solution.
  \end{theorem}

\noi{\em Proof.}
 It consists of a few steps.
  \underline{(i)  {\em Scaling invariance}}. One can see that \ef{20}
 admits a group of scalings: if $f_1(y)$ is a solution,
  \be
  \label{21a}
   \tex{
\mbox{then} \quad   f_a(y)= a f_1 \big( \frac y a\big) \quad
\mbox{is a solution for
   any $a \not = 0$}.
    }
    \ee

 \noi\underline{(ii) {\em Finite regularity at $y=0$}}.
 This can be seen by (still formal) expansion near the origin
  \be
  \label{22}
   f(y)=y + \e(y) \whereA \e(y)=o(y) \asA y \to 0.
    \ee
Substituting \ef{22} into the ODE \ef{20} and ``linearizing"
yields the following Euler's ODE:
 \be
 \label{23}
  \tex{
  y^2 \e''-2y\e'+\e=0 \LongA \e(y) =y^{m_\pm} \whereA m_\pm=\frac {3 \pm \sqrt
  5}2.
   }
   \ee
 Since we are looking for sufficiently smooth and at least
 classical solutions of the PDE \ef{15}, we have to choose the
 maximal exponent
  \be
  \label{23m}
   \tex{
   m=m_+=\frac {3 + \sqrt
  5}2= 2.61803... \in (2,3) \quad (m_-=0.38197...<1).
   }
   \ee
 The justification of the expansion \ef{22}, \ef{23} is performed
 in a standard manner by reducing the ODE to the equivalent
 integral equations for $y=y(f)$ (this deletes the trivial
 solution $f \equiv 0$) and applying Banach's Contraction
 Principle with a suitable weighted functional setting.

Thus, \ef{20} admits a 1D bundle of solutions close to the origin:
 \be
 \label{24}
 f(y) = y + A y^m+ o(y^m) \forA y >0,
  \ee
where $A \in \re$ is arbitrary. For $A=0$, we have the exact
unbounded solution
 \be
 \label{ff12}
 f_0(y)=y.
 \ee
 It follows from the symmetry \ef{19} that the extension of the
 shock profiles for $y<0$ is:
  \be
  \label{24-}
 f(y)=y+ A |y|^{m-1} y+o(|y|^m) \forA y<0.
 \ee

The parameter $A<0$ makes it possible to shoot the necessary
constant $+1$ at infinity, which is rather standard and we omit
details. $\qed$

\ssk

\begin{figure}
\centering
\includegraphics[scale=0.75]{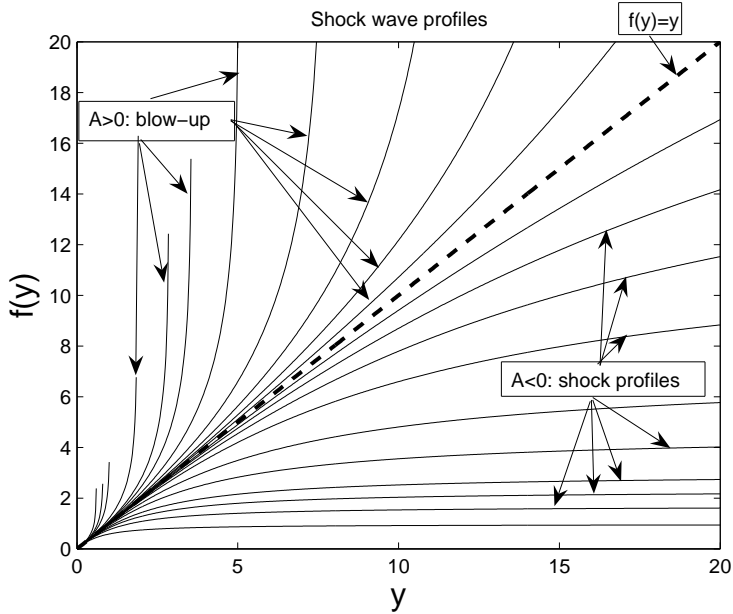}              
 \vskip -.3cm
\caption{\small  Various solutions of the ODE \ef{20} with the
expansion \ef{24} for different values of $A>0$.}
   \vskip -.3cm
\label{F1}
\end{figure}

Actually, some necessary properties of solutions of \ef{20} can be
easily checked numerically by using a standard {\tt MatLab}
solvers such as the {\tt ode45}, and this diminishes the necessity
of a deeper mathematical study. In Figure \ref{F1}, we show the
monotone dependence on $A<0$ and $A>0$ of solutions of \ef{20}
with the expansion \ef{24}. For $A>0$, the local solutions $f(y)$
blow-up in finite $y_0>0$ that is also easily proved. Note that
for $A<0$, according to \ef{152}, the behaviour in $\{y>0\}$ of
similarity solutions \ef{16} are locally governed by a well-posed
parabolic flow with a blow-up formation of a shock wave in finite
time.

The next Figure \ref{F2} demonstrates  a general view of
similarity shock profiles for both $y>0$ and $y<0$. By the
boldface line we denote the unique solution of the problem
\ef{20}, with the following values  in the expansion \ef{24}:
 \be
 \label{A11}
 A_*=-0.96...\,.
  \ee

\begin{figure}
\centering
\includegraphics[scale=0.75]{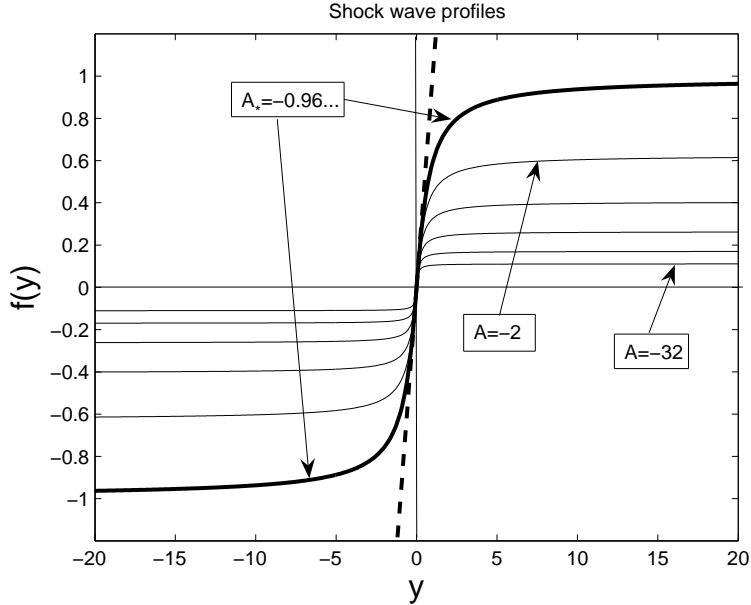}            
 \vskip -.3cm
\caption{\small  A general view of similarity shock profiles.}
   \vskip -.3cm
\label{F2}
\end{figure}

\begin{figure}
\centering
\includegraphics[scale=0.75]{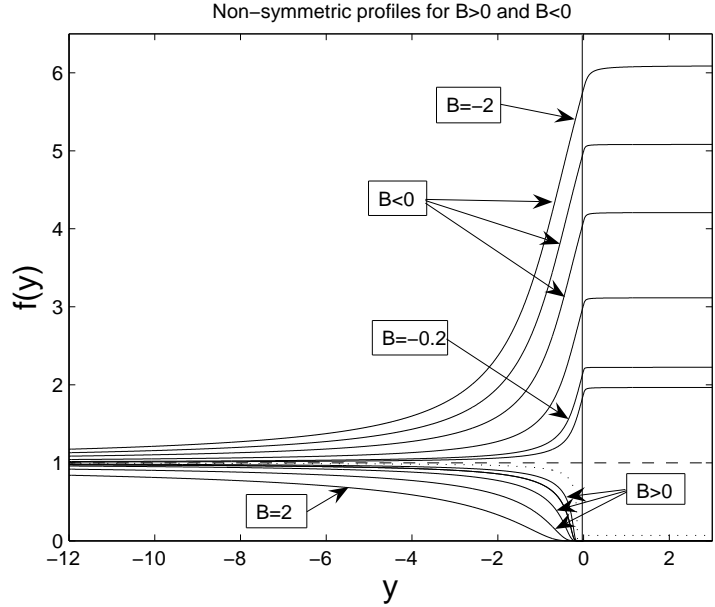}            
 \vskip -.3cm
\caption{\small  Examples of non-symmetric similarity shock
profiles satisfying \ef{bb5} for $B<0$ and $B>0$.}
   \vskip -.3cm
\label{F3}
\end{figure}

In addition, in Figure \ref{F3} we present the results of further
analysis of shock wave profiles as (generalized) solutions of the
ODE \ef{20}. Namely, we now shoot from $y=-50$ with the following
asymptotic expansion of the solutions about the equilibrium $f=1$:
 \be
 \label{bb5}
  \tex{
 f(y)=1 + \frac By + ... \asA y \to - \iy.
 }
 \ee
We now use the {\tt ode45} solver with the fully
 enhanced relative
and absolute tolerances
 \be
 \label{bb6}
 \mbox{Tols}=10^{-13}.
 \ee
For $B<0$, we observe that the solutions pass though the singular
point $y=0$ and converge to another larger equilibrium,
 \be
 \label{bb7}
 f(y) \to f_+ \asA y \to + \iy \whereA f_+=f_+(B)>1.
  \ee
Note that possible singularities at $y=0$ are integrable since
even in the worst case the solution is bounded,
 $$
 \tex{
 f'' \sim \frac 1y \LongA f(y) \sim C + y \ln |y|+... \asA y \to
 0,
 }
 $$
 or the regularity is better, when the right-hand side in \ef{20}
 vanishes at $y=0$ (this case actually occurs to give good solutions).

A different phenomenon is observed for $B>0$, where the shock
profiles seem have finite interfaces at some $y=y_0<0$. However,
next Figure \ref{F4} explains that this is wrong: thought (a) with
the enlargement $10^{-3}$ continues to convince existence of
negative finite interfaces, (b) given in the logarithmic scale
 confirms that the solutions have exponential decay up to $y=0^-$. Some lines in (b) are disjoint
 for profiles $|f(y)| \sim 10^{-17}$, which suddenly change sign due to numerical errors
 (recall the guaranteed accuracy \ef{bb6}).
 Namely, a more careful look at the equation \ef{20} shows that it
 admits solutions with the following typical non-analytic
 behaviour:
 \be
 \label{nn1}
  \tex{
   f(y) \sim {\mathrm e}^{-\frac a{(-y)}}\to 0 \asA y \to 0^-,
   }
   \ee
   where $a>0$ is arbitrary due to the scaling symmetry \ef{21a}.

Such solutions admit a natural trivial extension by $f(y) \equiv
0$ for $y \ge 0$, thus forming $C^\iy$ shock wave solutions
 with finite interfaces. Taking such a profile $f(y)$, instead of
 \ef{5}, we get in the limit the reflected  {\em Heaviside function}:
  \be
   \label{nn54}
   \lim_{t \to 0^-} u_{-}(x,t)= H(-x)= \left\{
    \begin{matrix}
     0 \quad \mbox{for} \quad x \ge 0,\\
  1 \quad \mbox{for} \quad x<0.
   \end{matrix}
    \right.
    \ee
It is principal to note that, due to the regularity \ef{nn1} at
the interface, the similarity solution \ef{16} is a $C^\iy$
function for $t<0$. In other words, \ef{nn54} describes formation
of a shock wave as $t \to 0^-$ from classic $C^\iy$ solutions.
 Such an evolution process of formation of shocks can be taken as
 a core property for introducing a test on the so-called {\em $\d$-entropy solutions},
 i.e., those which are obtained via smooth deformation of
 shocks appeared. See \cite{GPndeII, GalNDE5, GPnde}, where such an approach is
 developed for nonlinear dispersion equations such as
  \be
  \label{ND1}
  u_t=(uu_x)_{xx} \andA u_t=-(uu_x)_{xxxx}.
   \ee
Classic concepts of entropy do not apply to \ef{ND1}, so other
approaches are needed. On the other hand, a detailed pointwise
analysis of formation of gradient blow-up for PDEs such as
\ef{ND1} shows that uniqueness is violated there, so a consistent
entropy theory is most plausibly non-existent in principle; see
\cite{GalndeIII}. Deeper evolution and entropy-like
properties of such finite interface solutions are out
of the scope of
 this paper. Let us  mention that the questions on
 weak compactly supported solutions
 are in the focus of modern research in the EEs area; see
 \cite[p.~418]{Bar07}.


\begin{figure}
\centering
\subfigure[Enlargement $10^{-3}$]{
\includegraphics[scale=0.52]{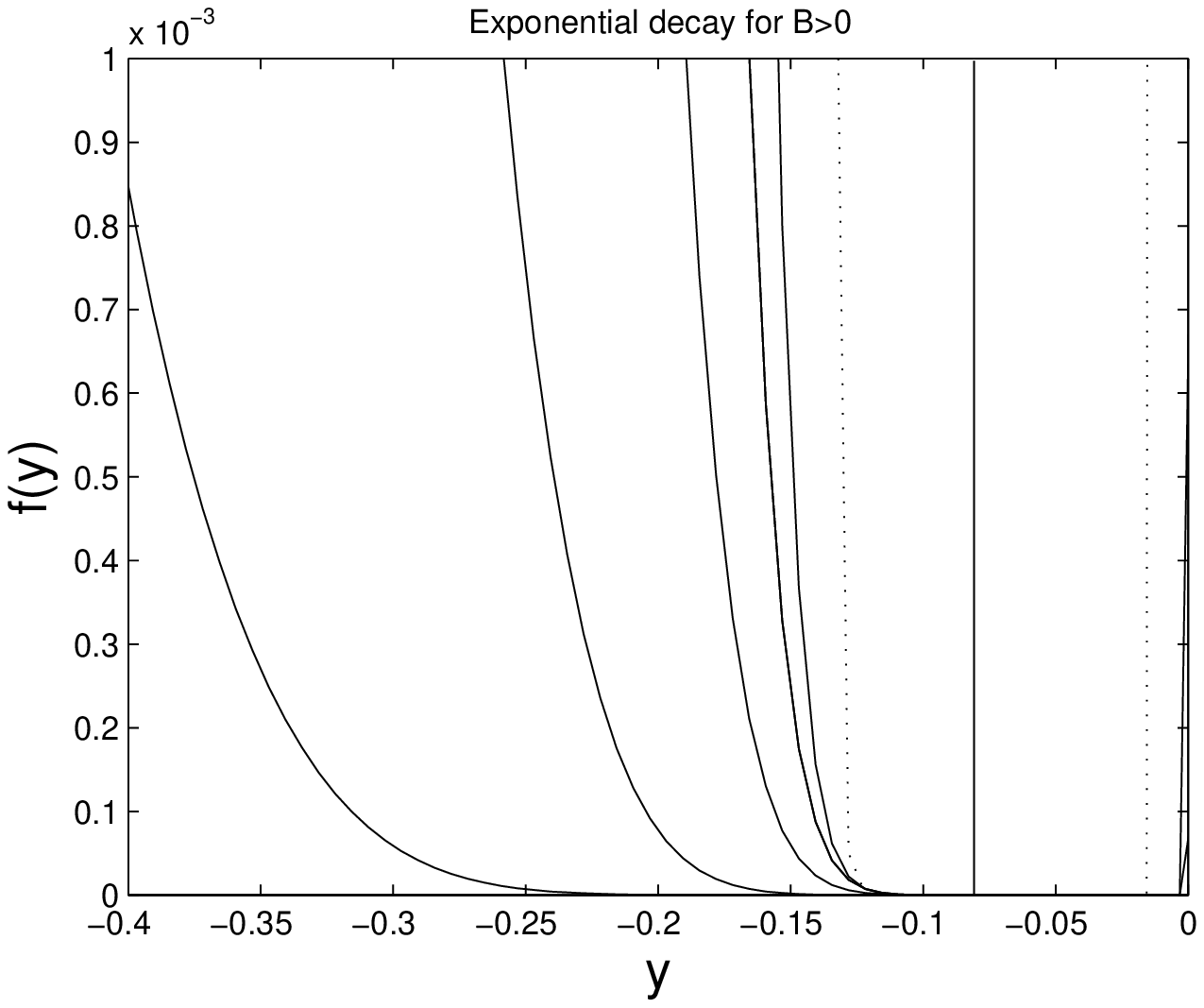}
}
\subfigure[Logarithmic scale]{
\includegraphics[scale=0.52]{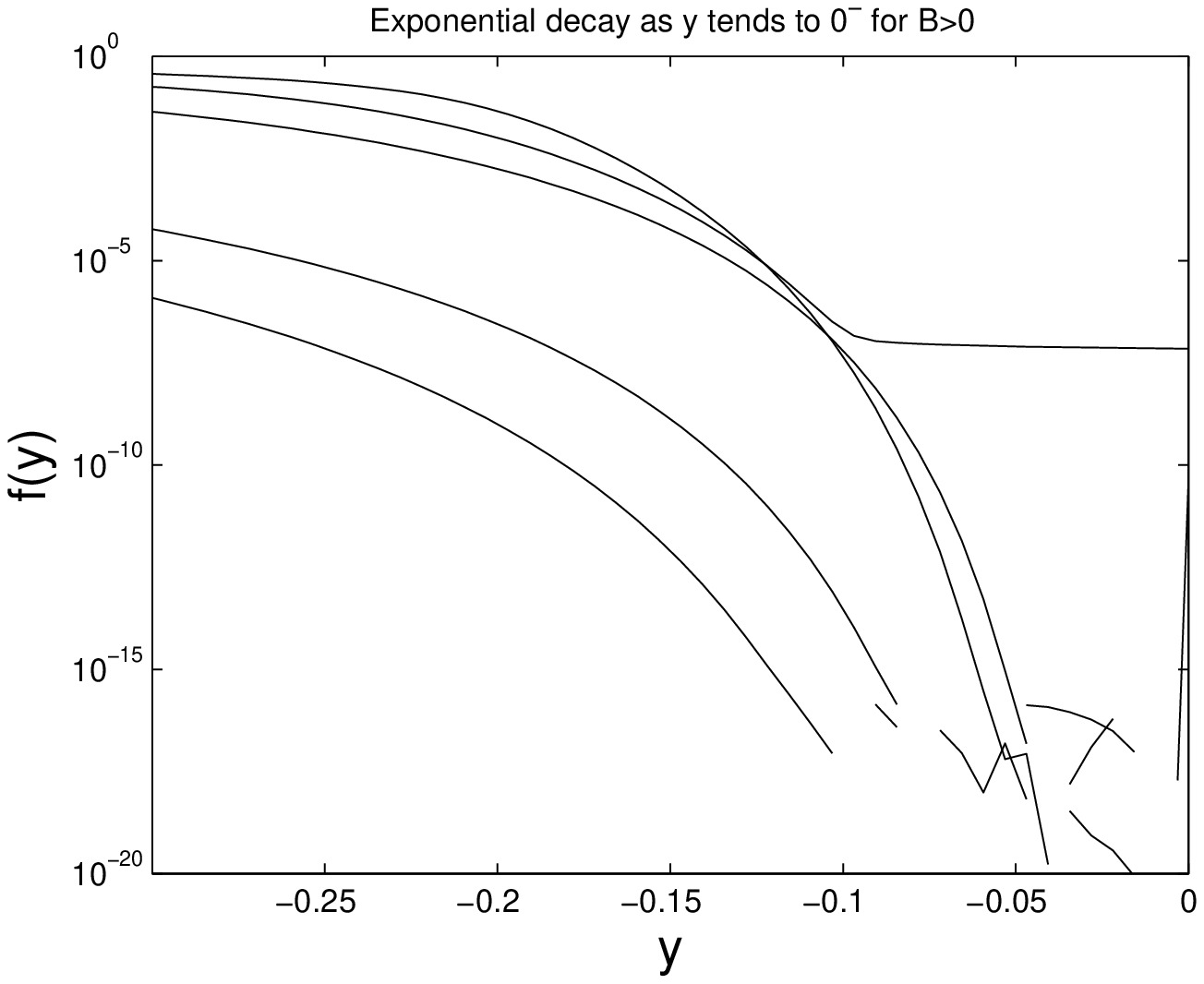}
}
 \vskip -.4cm
\caption{\rm\small  Enlarged shock profiles for $B>0$: exponential
decay as $y \to 0^-$.}
 \vskip -.3cm
 \label{F4}
\end{figure}


\subsection{Existence of rarefaction waves}

The corresponding rarefaction waves describe collapse of initial
singularities. The similarity ones then are defined for all $t>0$:
\be
 \label{16R}
  \tex{
  u_{+}(x,t)= F(y), \quad y= \frac x{t}, \quad
   y F F''=(F')^2(y-F)-2 F F', \quad F(\pm \iy)= \mp 1.
   }
   \ee
 It follows that, in comparison with the shock profiles $F(y)$
 constructed above, the rarefaction ones are given by the
 reflection:
 \be
 \label{81}
 -F(y) \quad \mbox{is a solution of the rarefaction ODE (\ref{16R}).}
 \ee
Therefore, the similarity solution in \ef{16R} is formed by the
initial singularity $S_-(x)=-{\rm sign}\, x$ and describes its
evolution collapse for $t>0$.

 Similarly, for data $u(x,0)=-H(-x)$, we obtain  similarity
 solution \ef{16R} describing its collapse. To get a unique
 solution, two initial functions should be prescribed.
 It is easy to see that $u_t(x,0)$ is not a distribution.
 Indeed, we have,
  according to \ef{16R} and asymptotics \ef{bb5}, \ef{nn1}, in the sense of
 distributions, it is formally valid that
  \be
  \label{ff2}
   \tex{
  u_t(x,t)= - \frac 1t \, F'(y) y \to a_0 \d(x) \asA t \to 0^+
  \whereA a_0= - \int F'(y) y\,{\mathrm d}y > 0.
 }
   \ee
 However, in view of the asymptotics \ef{bb5}, the integral in
 \ef{ff2} diverges, so that $u_t(x,0)$ is not a bounded positive
 measure. In the sense of distributions, such a behaviour can be
 characterized by a logarithmic divergence as follows:
  \be
  \label{ff3}
   \tex{
  u_t(x,t)= - \frac 1t \, F'(y) y \sim B |\ln t|\, \d(x) \asA t
  \to 0^+,
  \quad\mbox{so \,$\frac {u_t}{|\ln t|}(x,0)$\, is a measure}.
  }
   \ee

\section{Unbounded shocks by similarity blow-up and rarefaction profiles}
\label{S2N}

More complicated shock-type structure occurs if in \ef{16} we
introduce a parameter $\a \in (0,1)$ looking for similarity
solutions (cf. \ef{sw1})
 \be
 \label{16N}
   \tex{
  u_{-}(x,t)=(-t)^{-\a} f(y), \quad y= \frac x{(-t)^\b} \forA t<0
   \whereA
  \b=1-\a,
   }
   \ee
 and $f(y)$ then solves the following more complicated ODE:
 \be
 \label{17N}
 \b^2 y^2 f f''=\b y(f')^2(\b y+f)+ \b(\a-2)y f f'+ \a f^2(f'-1)
 \inB \re.
   \ee
   Using the same symmetry \ef{19},
 we again pose the anti-symmetry condition at the origin
  \be
  \label{as1}
  f(0)=0 \LongA f(-y) \equiv -f(y).
   \ee
In view of  \ef{13}, we have to require that
 \be
 \label{ss1A}
  \tex{
  \frac{u_t}u \equiv (-t)^{-\a-1} \frac{\a f+ \b f' y}f>0,
  }
  \ee
  which for \ef{17N} can be also checked by a Maximum Principle
  approach. One can see that this means that $u_t>0$ for $x>0$
  implying a monotone in time growth of blow-up solutions
  (a typical feature of parabolic equations, see \ef{152}).
 By a similar shooting technique, we prove the following analogy
 of Theorem \ref{T1}, which includes all necessary asymptotics of blow-up profiles:

 \begin{theorem}
 \label{T2}

 The problem $\ef{17N}$, $\ef{as1}$ admits a family of solutions
 $f(y)>0$ for $y>0$ satisfying for any $A<0$:
  \be
  \label{as2}
   \tex{
  f(y)=y + Ay^{m_+}+ ... \asA y \to 0^+, \quad m_\pm=
   \frac{2\a^2-4\a+3 \pm
   \sqrt{(2\a^2-4\a+3)^2-4(1-\a)^4}}{2(1-\a)^2};
 }
    \ee
 \be
 \label{as3}
 f(y) = C_0 y^{-\frac \a{1-\a}}+... \to 0^+ \asA y \to + \iy \whereA
 C_0=C_0(A)>0.
 \ee
 \end{theorem}

Let us note that the linearization \ef{22} now leads to the
equation
 \be
 \label{23N}
  \tex{
  (1-\a)^2 y^2 \e'' - (\a^2-2\a+2)y \e'+(1-\a)^2\e=0,
  }
 \ee
whence the same solutions with the exponents $m_\pm$ given in
\ef{as2}.

It follows that \ef{16N} is an {\em unbounded} blow-up solution
that forms as $t \to 0^-$ the following singular {\em final time
profile} (cf. \ef{18}):
 \be
 \label{as4}
 u_{-}(x,0^-)= \left\{
 \begin{matrix}
 \,\,\,
C_0 x^{-\frac \a{1-\a}}\,\, \,\,\,\forA x>0,\\
   -C_0 |x|^{-\frac \a{1-\a}} \forA x<0
   \end{matrix}
   \right.
   \quad ( u_{-}(0^\pm,0^-)=\pm \iy).
    \ee
The asymptotics at infinity \ef{as3} is also a standard one, which
is approved as in reaction-diffusion theory; see e.g.
\cite[Ch.~4]{SGKM} and references therein.

\begin{figure}
\centering
\includegraphics[scale=0.75]{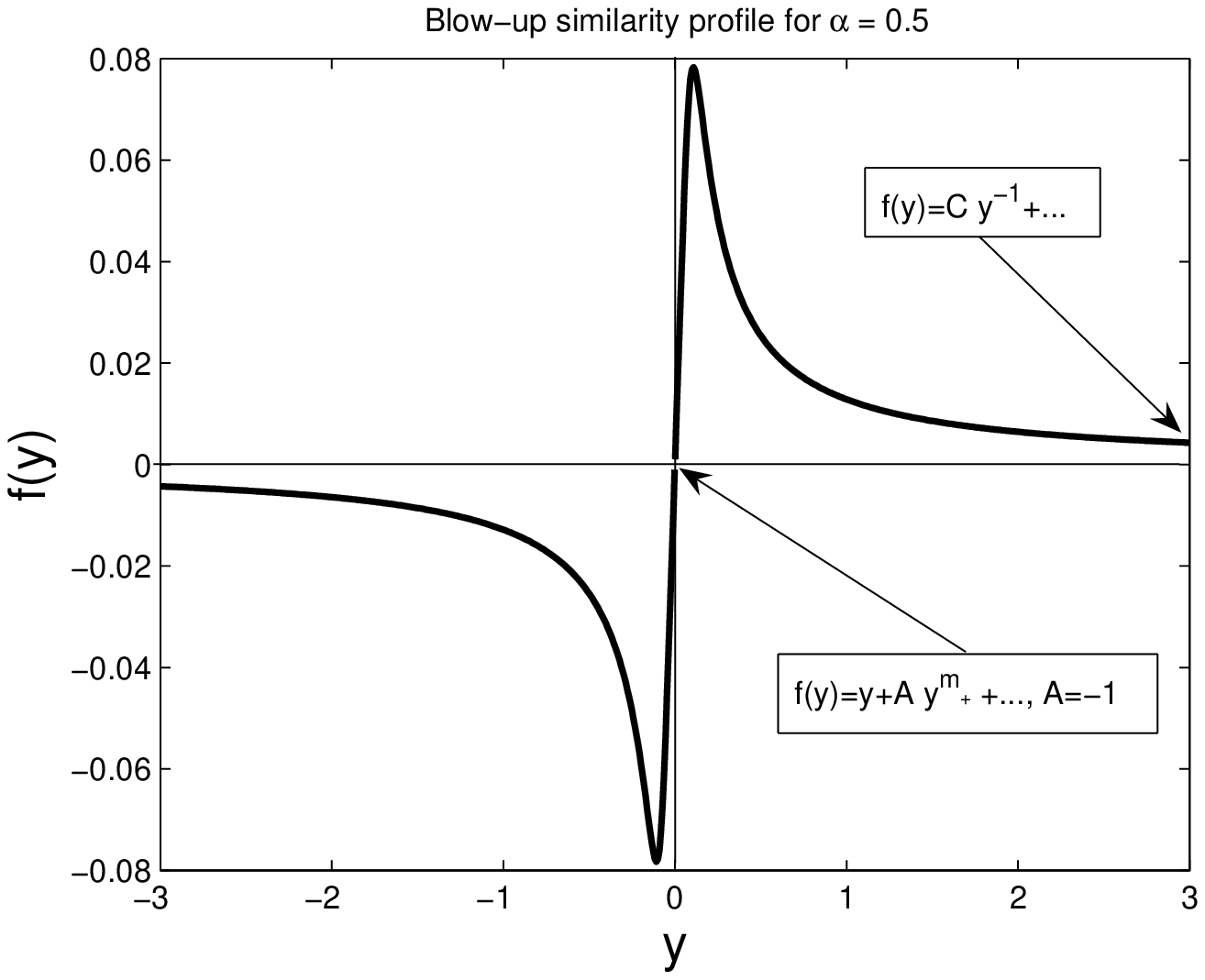}              
 \vskip -.3cm
\caption{\small  A solutions of  \ef{17N}, \ef{as1} for $\a= \frac
12$ and  $A=-1$.}
   \vskip -.3cm
\label{F1NN}
\end{figure}
\begin{figure}
\centering
\includegraphics[scale=0.75]{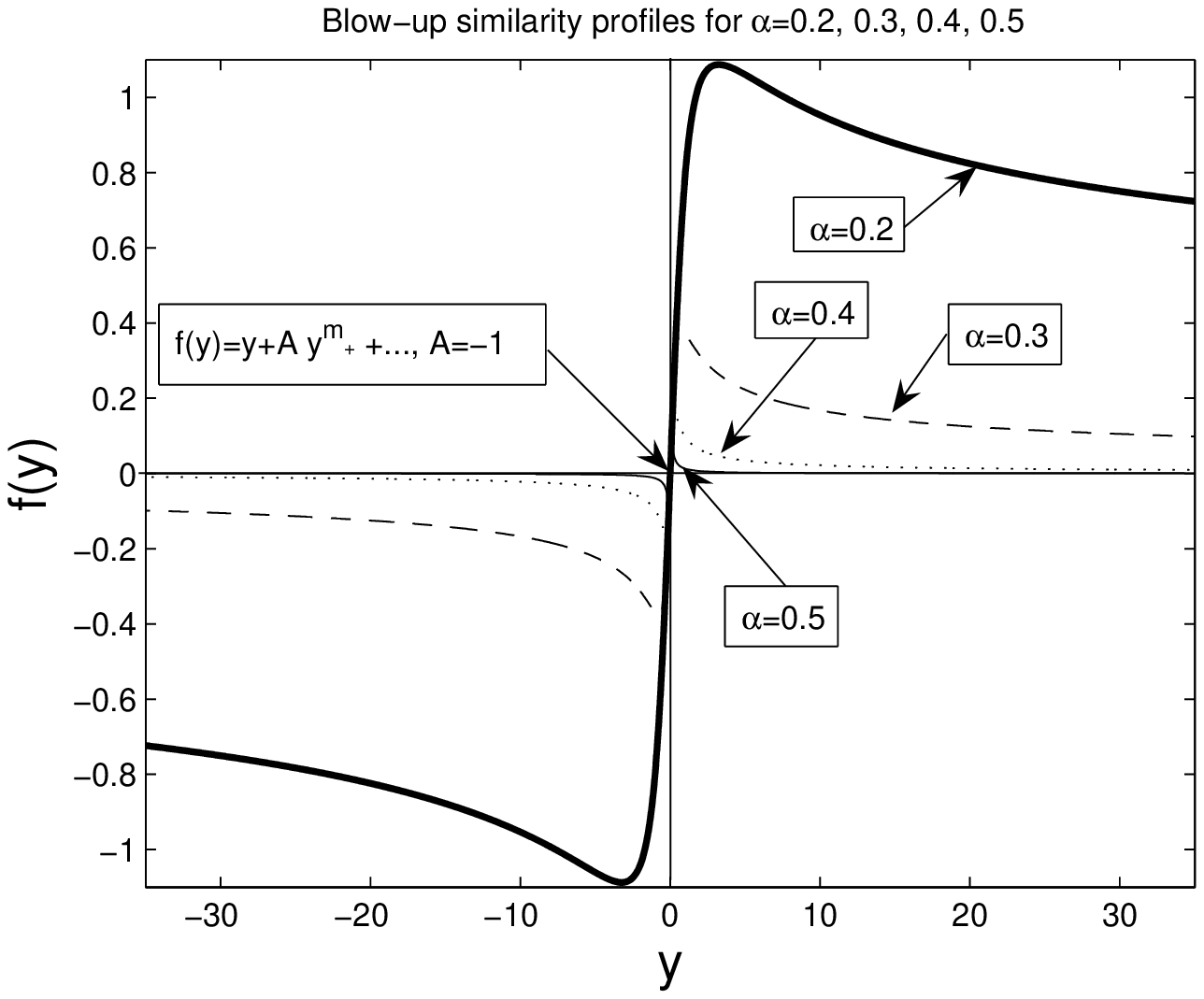}              
 \vskip -.3cm
\caption{\small  Solutions of  \ef{17N}, \ef{as1} for $\a=0.2$,
$0.3$, $0.4$, and $0.5$;  $A=-1$.}
   \vskip -.3cm
\label{F2NN}
\end{figure}

As usual, we complete our analysis by numerics. In Figure
\ref{F1NN}, we show a typical blow-up profile for $\a= \frac 12$,
while Figure \ref{F2NN} explains deformation of such profiles
$f(y)$ with $\alpha$.

\subsection{Rarefaction waves from singular initial data}

Similar to \ef{16R}, the global similarity solutions of \ef{15},
 \be
 \label{16NR}
   \begin{matrix}
  u_{+}(x,t)=t^{-\a} F(y), \quad y= \frac x{t^\b} \forA t>0,
  \quad
  \b=1-\a,\qquad\quad \ssk\ssk\\
  \b^2 y^2 F F''=\b y(F')^2(\b y-F)+ \b(\a-2)y F F'+ \a F^2(-F'-1)
 \inB \re, \qquad\quad
   \end{matrix}
   \ee
 corresponding to the reflection $F \mapsto -F$ in the blow-up ODE
 \ef{17N}, describe collapse of singular initial data \ef{as4}
 posed at $t=0$. We then again arrive at a difficult question on
 the correct entropy-like choice of proper solutions, which reveal
  some remnants of those ones for the EEs \ef{NS1} and hyperbolic systems
 such as \ef{HS1}. Hopefully, an evolution entropy test then can
 be developed in lines of smooth $\d$-deformations as in
 \cite{GPndeII, GalNDE5, GPnde}. However, there are examples of
 higher-order PDEs, for which any entropy-like unique extension of
 a solution after singularity is principally impossible (there
 are infinitely many extensions that are all equally allowed),
 \cite{GalndeIII}. We do not know whether such negative trends
 can be oriented to the EEs \ef{NS1}.

 \section{Uniqueness of a local continuation after gradient blow-up}
  \label{S31}

We now need to consider the principal question on a generic
(self-similar) mechanism of formation of shocks for the model
\ef{15}. To this end, following \cite{GalndeIII}, we will use the
same similarity solutions \ef{16N}, but now with
 \be
 \label{al1}
  \tex{
 \a<0 \LongA \b=1-\a \equiv 1+|\a|>1, \quad \mbox{so that}
 \quad \frac {|\a|}\b \equiv \frac {|\a|}{1+|\a|} <1.
 }
 \ee
 It can be shown that such solutions \ef{16N}  create
 continuous ``initial data"  (cf. \ef{as4})
 \be
 \label{as4N}
 u_{-}(x,0^-)= \left\{
 \begin{matrix}
 \,\,\,
C_0 x^{\frac {|\a|}{1+|\a|}}\,\, \,\,\,\forA x>0,\\
   -C_0 |x|^{\frac {|\a|}{1+|\a|}} \forA x<0,
   \end{matrix}
   \right.
    \ee
 where $C_0>0$ is some (actually arbitrary by scaling \ef{21a}) constant.
 Note that \ef{as4N} is still continuous at $x=0$ but has a first
 gradient blow-up there.

The question is how to extend the data \ef{as4N} for $t>0$. This
is naturally done by using the global similarity solution
\ef{16NR}. Unlike the previous cases, we here assume that a shock
occurs at $y=0$, i.e., we put two conditions there
 \be
 \label{co1}
 F(0)=F_0>0 \andA F'(0)=F_1.
  \ee
By a rather standard local analysis of the ODE \ef{16NR} at $y=0$,
it is not difficult to conclude that the Cauchy problem \ef{co1}
has a unique smooth solution in the only case
 \be
 \label{co2}
  \tex{
  F'(0)=F_1=-1 \LongA F(y)=F_0-y+ \frac {\b^2 y^2}{2 \a F_0}+...\,
  .
   }
   \ee
For other values of $F'(0) \ne -1$, the solution is not smooth and
is singular at $y=0$. All such singular bundles can be also
derived. Therefore, the behaviour of solution $F(y)$ as $y \to
0^+$ is very unstable and, even numerically, its construction is
not straightforward.

\begin{figure}
\centering
\includegraphics[scale=0.75]{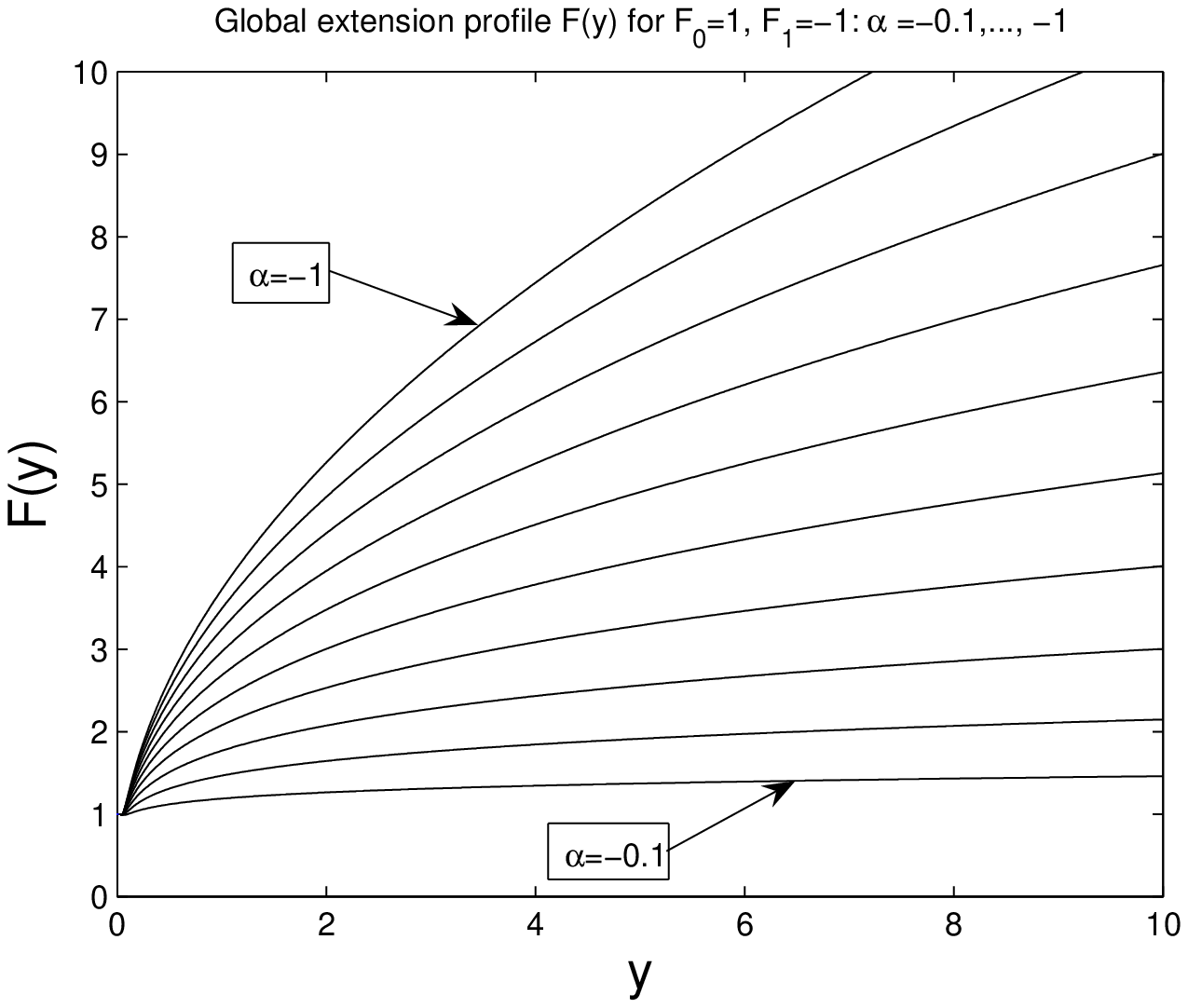}              
 \vskip -.3cm
\caption{\small  Solutions of the ODE \ef{16NR}, \ef{co2} for
$\a=-0.1$, $-0.2$,... $-1$.}
   \vskip -.3cm
\label{FGl1}
\end{figure}

 Before stating the main result, we present in Figure \ref{FGl1}
such global similarity profiles $F(y)$ for various values of $\a
\in [-1,-0.1]$. Since the non-monotonicity of $F(y)$ is not
visible here at all, in the next Figure \ref{FGl2}, we show
typical shape of such profiles in the Cauchy problem posed for $y
>\d=0.1$, with $F(\d)=1$ and various $F'(\d)<0$. Note that the
absolute minimum points at some $y=y_0>0$ perfectly corresponds to
the Maximum Principle for the ODE \ef{16NR}. Indeed, one has:
 \be
  \label{co3}
  F'(y_0)=0, \,\, F(y_0)>0 \LongA \b^2 y_0^2
  F(y_0)F''(y_0)=|\a|F^2(y_0)>0.
   \ee

\begin{figure}
\centering
\includegraphics[scale=0.75]{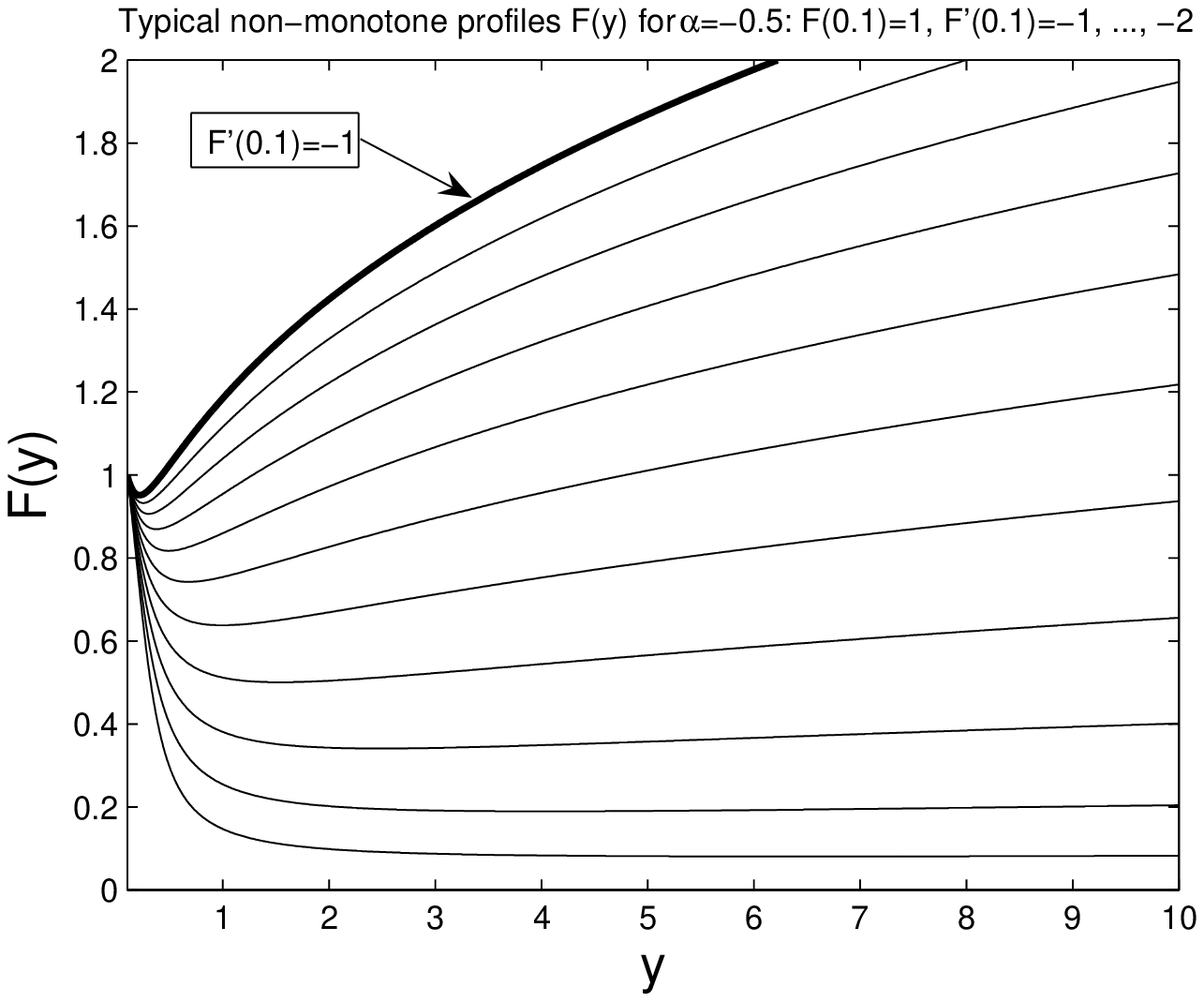}              
 \vskip -.3cm
\caption{\small  Typical non-monotone solutions of the ODE
\ef{16NR}, \ef{co2} for $\a=-0.5$.}
   \vskip -.3cm
\label{FGl2}
\end{figure}

By checking the asymptotic properties of solutions of \ef{16NR} as
$y \to + \iy$, it is not difficult to see that all the profiles
with the expansion \ef{co2}, which do not vanish by the Maximum
Principle shown in \ef{co3}, exhibit the desired behaviour at
infinity:
 \be
 \label{co4}
F(y)=C y^{\frac {|\a|}{1+|\a|}}(1+o(1)) \asA y \to +\iy.
 \ee
Here $C$ is a positive constant but not necessarily the one $C_0$
obtained in \ef{as4N} by the blow-up limit. However, the
transition $C \mapsto C_0$ is uniquely done by the scaling
invariance \ef{21a}, so we arrive at:

\begin{proposition}
 \label{Pr.Ext}
 For given initial data $\ef{as4N}$, there exists a unique
 similarity profile $F(y)$ such that the global similarity
 solution $\ef{16NR}$ satisfies
  \be
  \label{co5}
  u_{-}(x,0^-)=u_{+}(x,0^+) \inB \re.
   \ee
   \end{proposition}

We call the corresponding profiles $\{f,F\}$ a
 {\em similarity
global extension pair}. Thus, this pair is always unique, so the
self-similar  blow-up solution $u_-(x,t)$ for $t<0$ has a unique
extension via the similarity one $u_+(x,t)$ for $t>0$.

This justifies a kind of ``micro-local" uniqueness theory of shock
wave solutions for the PDE \ef{15}, which sounds rather
optimistic. Existence, or even a possibility of construction, of a
global ``entropy-like" theory for \ef{15} (say, in
Oleinik--Kruzhkov--Lax--Glimm--Bressan--etc. sense) is quite
questionable still. In any case, the  above ``micro-uniqueness"
result is rather inspiring bearing in mind that this could affect
some similar extension properties of the EEs.

\section{Final conclusions: towards the full model of $(r,\var)$-shock and blow-up formation}
\label{S4}

The integral equation \ef{11} assumes another mechanism of shock
formation for data:
 \be
 \label{b1}
V_{0\var} \approx 0 \andA V_0 \approx 1.
 \ee
 Then the second term dominates and, on two differentiations, this
 leads to the PDE:
  \be
  \label{b2}
   U_{r\var}= \Big[ \frac{U_{r\var}}{\frac U{U_t}\big( \frac
   {U_t}U\big)_t- U_r}\Big]_t.
   \ee
 This equation is supposed to describe more complicated phenomena
 of shock and more general blow-up formation, where both independent variables $(r,\var)$
 are essentially and equally involved. Of course, \ef{b2} is
 uncomparably more difficult than \ef{15}, and seems can give a deeper insight
 into weak and strong singularity formation for the EEs (and also for \ef{HS1}) in the actual
  2D $\{r,\var\}$-geometry.
  We do not study this equation
 here. Nevertheless, it should be mentioned that equations such as
 \ef{b2} and others must be carefully investigated
  and understood before any
 serious attack of singularity phenomena for Euler equations in 2D
 and, especially,  3D, which will indeed require other hierarchy of some reduced PDE models to
 introduce and study.





\end{document}